# New artificial-free phase 1 simplex method

Nasiruddin Khan, Syed Inayatullah*, Muhammad Imtiaz

and Fozia Hanif Khan

Department of Mathematics, University of Karachi,

Karachi, Pakistan, 75270

*Email: inayat@uok.edu.pk

**Abstract:**

*This paper presents new and easy to use versions of primal and dual phase 1 processes which obviate the role of artificial variables and constraints by allowing negative variables into the basis. During the process new method visits the same sequence of corner points as the traditional phase 1 does. The new method is artificial free so, it also avoids stalling and saves degenerate pivots in many cases of linear programming problems.*

**\*Corresponding Author**

## 1. Introduction:

Initial requirement of the simplex method is a basic feasible solution and whenever an initial basic feasible solution of an LP is not given, we should apply the simplex method in two phases [4,10], called phase 1 and phase 2. In phase 1 we create a basic feasible solution artificially by adding some (non-negative) artificial variables to the problem with an additional objective, equal to minimization of the sum of all the artificial variables, called infeasibility form. Here in this paper we call it phase 1 objective. The purpose of phase 1 process is to maintain the feasibility and minimize the sum of artificial variables as possible. If phase 1 ends with an objective value equal to zero, it implies that all artificials have been reached to value zero and our current basis is feasible to the original problem, then we may turn to the original objective and proceed with simplex phase 2. Otherwise we conclude that the problem has no solution.

The above approach is the most traditional but not the only one to do so. In [11] Zoutendijk presents different variations of the phase 1 simplex method. On page 47 he has also presented an artificial free Big M like method. Recently Arsham [1,2] proposed alternate but artificial-free methods to perform phase 1. In section 6 of the present working we shell present another artificial-free version of phase 1 process which is a modified form of the method presented in [11] and hence quite equivalent to the traditional phase 1 process of artificial variables  (described above) because it visits the same sequence of corner points as the traditional phase 1 does. But the advantage of new approach is that it could start with an infeasible basic solution, so introduction of artificial variables is not mandatory. Additionally, our new phase 1 has another major advantage over the traditional approach that is the traditional phase 1 may encounter problem of stalling due to degenerate artificial variables but new method resolves the problem effectively and saves those degenerate pivots. In section 10 we shell also describe the dual counterpart of our new artificial free phase 1 which is indeed an artificial constraint free version of traditional dual simplex phase 1.

## 2. The Basic Notations:

A general LP problem is a collection of one or more linear inequations and a linear objective function involving the same set of variables, say $x_1, \ldots, x_n$, such that

$$\text{Maximize} \quad c^T x$$
$$\text{subject to} \quad Ax = b \,,$$
$$x \geq 0, \, x \in \Re^n$$

where $x$ is the decision variable vector, $A \in \Re^{m \times n}$, $b \in \Re^m$, $c \in \Re^n$.

We shell use $a_i$ and $a_{.j}$ to denote the $i^{\text{th}}$ row and $j^{\text{th}}$ column vector of $A$. Now we define a basis $B$ as an index set with the property that $B \subseteq \{1,2,\ldots n\}$, $\mid B \mid = m$ and $A_B := [a_{.j} \mid j \in B]$ is an invertible matrix, and Non-basis $N$ is the complement of $B$. i.e. $N := \{1,2,\ldots n\} \backslash B$.

So, we may construct the following dictionary for basis $B$ c.f. [3],

$$D(B) = \begin{bmatrix} z & -\bar{c}^T \\ \hline \bar{b} & \bar{A} \end{bmatrix}$$

where,

$$\bar{A} = A_B^{-1} A_N \,,$$

$$\bar{c}^T = c_N^T - c_B^T A_B^{-1} A_N \,,$$

$$\bar{b} = A_B^{-1} b$$

$$z = C_B^T A_B^{-1} b$$

The associated basic solution could directly be obtained by setting $x_B = \overline{b}$. Here onward in this text we assume that the reader is already familiar about pivot operations, duality and primal-dual relationships. It is well known that if $d_{B0} \geq 0$ then $B$ (or $D(B)$) is called primal feasible (dual optimal); if $d_{0N} \geq 0$ then $B$ (or $D(B)$) is called primal optimal (dual feasible); and if both $d_{B0} \geq 0$ and $d_{0N} \geq 0$ then $B$ (or $D(B)$) is called optimal feasible. A basis $B$ (or a dictionary) is called inconsistent if there exists $i \in B$ such that $d_{i0} < 0$ and $d_{iN} \geq 0$, and unbounded if there exists $j \in N$ such that $d_{0j} < 0$ and $d_{Bj} \leq 0$.

### 3. Formation of Auxiliary form for traditional simplex method:

Through out the paper we shell call the following form of LP problem as standard form, if $A \in \Re^{m \times p}$, $b \in \Re^{m}$, $c \in \Re^{p}$.

Maximize $\quad c^{T}x$

subject to $\quad Ax \leq b$ ,

$\quad\quad\quad\quad x \geq 0, x \in \Re^{p}$.

Here $b$ would not necessarily be completely non-negative. By adding the slack vector $s$, we can have an equivalent equality form of the above system,

Maximize $\quad c^{T}x$

subject to $\quad A x + s = b$,

$\quad\quad\quad\quad x \geq 0, s \geq 0$.

$\quad\quad\quad\quad x \in \Re^{p}, s \in \Re^{m}$.

Let $S$ is index set of variables in $s$. Clearly, for above system the readily available basis is $S$. But $S$ may not constitute an initial feasible basis for simplex method. To follow the traditional simplex phase 1 approach we must transform such system into the following auxiliary form,

Minimize     $1v$

Maximize     $c^T x$

subject to    $A x + s - v = b$

$x \geq 0,\, s \geq 0,\, v \geq 0$

$x \in \Re^p,\, s \in \Re^m,\, v \in \Re^m$

 Purpose of the new objective function *Minimize 1v* (called the phase 1 objective or

infeasibility form) is to force the artificial vector $v$ to zero, because $Ax + s = b$ if and

only if $Ax+s-v=b$ with $v =0$. It is important to note that for each slack variable in $s$

there is an artificial variable (of opposite sign) in $v$, so we may construct a one to one

correspondence of each variable in $s$ with a negative conjugate in $v$. In the above

system where ever feasible we would take variables from the slack vector $s$ as basic

variable and remaining basis would be formed from $v$.  The initial non-basic variables

of $v$ are the permanent non-basic variables and the remaining are temporarily basic

variables but once they leave they would become permanent non-basic. Clearly in the

first phase, the above auxiliary system either provides a feasible basis or shows that

the original system has no feasible basis. For proofs of correctness and details see

[12].

**4. A useful trick to reduce the computational efforts due to degenerate pivots in**

**Simplex Phase 1:**

Whenever the dictionary encounters the degeneracy in artificial variables there is a way to drive the degenerate artificial variables out of the basis by placing a legitimated variable from $x$ or $s$ into the basis. As stated in the previous section, for each basic artificial variable there is a slack variable with coefficient -1, which is hence eligible to enter the basis.  Now it is straight forward to say that, from the efficiency point of view the best choice of placement of artificial variables is by that associated slack variable. The reason is quite clear because by performing this type of degenerate pivot the whole dictionary, except the pivot row and objective row, would remain preserved. More over the pivot row would only need to be just multiplied by '-1'.

**Lemma:** One could reduce the computational effort in degenerate pivots, due to degeneracy in artificial variables, during simplex method by making the artificial variable as leaving and corresponding slack as entering basic variable.

## 5. Deduction of artificial free form from Auxiliary form of LP.

We decompose $s$ and $v$ into $s_{M_1}$, $s_{M_2}$, $v_{M_1}$, $v_{M_2}$, where index sets $M_1$ and $M_2$ are chosen such that $b_{M_1} \geq 0$ and $b_{M_2} < 0$, respectively. We can have an equivalent to Auxiliary form mentioned in section 3,

$$\text{Minimize} \quad \sum_{i \in M_2}(a_i x + s_i - b_i)$$

$$\text{Maximize} \quad c^T x$$

$$\text{Subject to} \quad A\,x + s - v = b,$$

$$x \geq 0,\ s_{M_1} \geq 0,\ s_{M_2} = 0,\ v_{M_1} = 0,\ v_{M_2} \geq 0.$$

$$x \in \Re^p,\ s \in \Re^m,\ v \in \Re^m.$$

It is quite clear that vector $v_{M_1}$ are permanent non-basic variables and the remaining $v_{M_2}$ are currently in the basis but once they leave the basis would never come in again. We can introduce a new variable vector $w = s - v$ ,which has the property that when ever $v$ is positive it would be negative, that is showing infeasibility of the current basis.

$$\text{Minimize} \quad \sum_{i \in M_2} (a_i x + s_i - b_i)$$

$$\text{Maximize} \quad c^T x$$

subject to     $A\,x + w = \ b,$

$x \geq 0,\ s \geq 0,\ \ w$ is unrestricted.

$x \in \Re^p,\ s \in \Re^m, w \in \Re^m.$

It should be noted that here in phase 1 objective the coefficient of slack variables $s_{M_2}$ are always (unit) positive. So during phase 1 any variable of $s_{M_2}$ would not be a candidate for entering basic variable. More over the quantity $1 b_{M_2}$ is insignificant to us. Hence we can safely remove the terms $1 s_{M_2}$ and $1 b_{M_2}$ form the phase 1 objective function. And the new artificial free Auxiliary form of the given LP is,

$$\text{Minimize} \quad \sum_{i \in M_2} a_i x$$

$$\text{Maximize} \quad c^T x$$

subject to     $A\,x + w = \ b,$

$x \geq 0,\ \ w$ is unrestricted.

$x \in \Re^p,\ w \in \Re^m.$

## 6. The new artificial free phase 1 :

### Problem 1:

Given a dictionary $D(B)$, obtain primal feasibility.

### Algorithm :

**Step 1:** Let $L$ be a maximal subset of $B$ such that $L = \{i : d_{i0} < 0 , i \in B\}$. If $L = \emptyset$ then

    $D(B)$ is primal feasible. **Exit.**

**Step 2:** Compute phase 1 objective vector $W(B) \in \Re^N$ such that $W(B)_j = \sum_{l \in L} d_{lj}$ .

**Step 3:** Let $K \subseteq N$ such that $K = \{j : W(B)_j < 0, j \in N\}$. If $K = \emptyset$ then $D(B)$ is primal

    infeasible. **Exit.**

**Step 4:** Choose $m \in K$ such that $W(B)_m \leq W(B)_k \ \forall \ k \in K$

    (Ties should be broken arbitrarily)

**Step 5:** Choose $r \in B$ such that

    $r = \text{argmin} \ \{d_{i0}/d_{im} \mid (d_{i0} < 0 , d_{im} < 0) \text{ or } (d_{i0} \geq 0 , d_{im} > 0), i \in B \}$

**Step 6:** Make a pivot on $(r,m)$. ($\Rightarrow$ Set $B := B + \{m\} - \{r\}$, $N := N - \{m\} + \{r\}$ and

    update $D(B)$.)

**Step 7:** Go to **Step 1**.

## 7. Explanation:

The basic strategy of our approach is to increase the number of feasible basic variables subject to preserve the feasibility of the existing feasible variables. For $m \in N, \ r \in B$, the entering basic variable $x_m$ and pivot column $d_{Bm}$ are could be determined by applying Dantzig's largest coefficient rule [5], Steepest edge pivot rules [6,7,8], largest distance pivot rule [9] or any other appropriate pricing rule on $S$.

Throughout the paper we have chosen largest coefficient pivot rule for entering basic variable. As soon as the non-basic variable $x_m$ becomes basic some of the existing basic variables would *increase* and some of them would *decrease*. So, we may divide the current basic variables into four categories *infeasible & increasing, infeasible & decreasing, feasible & increasing, and feasible & decreasing*. The leaving basic variable $x_r$ would be the basic variable which firstly increases or decreases to zero. Clearly $x_r$ is possible to leave only when it is either infeasible & increasing or feasible & decreasing. This kind of leaving basic variable could be determined by simply taking minimum ratio test. Here the procedure of minimum ratio test (see [11]) is different from the traditional method. In this process as described in **Step 5** we take ratios of right hand side of feasible constraints with corresponding element in the pivot column only when the denominator is a positive element and for infeasible constraints only when the denominator is a negative element.

After determining entering and leaving basic variables next step is to update the basis and the associated dictionary. Just like traditional phase 1 the new version will continue until all constraints become feasible, and then if needed we may turn to usual phase 2 process to reach the optimality.

## 8. Proof of Correctness:

Our artificial free phase 1 could start with an infeasible basis without making it artificially feasible. As stated earlier, at the end of each iteration objective of the traditional phase 1 is to minimize of the sum of all artificial variables remained left in the basis. In an equivalent sense, as shown in section 6, at the end of each iteration we compute the phase 1 objective vector $W$ of our artificial free approach explicitly by computing sum of infeasible constraints (constraint with negative right hand value). And just like traditional simplex, our method intends to achieve the feasibility of the infeasible variables subject to preserve the feasibility of existing feasible variables.

It is easy to realize that our method is just a simplified image of the traditional method. The only difference occurs in degenerate pivots, due to degenerate artificial variables, where the new method skips that pivot. If for instance we make a relationship between degenerate variables of both the traditional and the artificial free dictionaries, we may conclude that each degenerate leaving artificial variable corresponds to a degenerate increasing variable in artificial free dictionary. In the traditional phase 1 process we may have to perform a degenerate pivot to make that degenerate artificial variable out of the basis, but as described in **step 5** in the our artificial free approach we do not allow this kind of pivot and should look for next minimum ratio. The reason is quite clear and well justified, because in such a particular case of degeneracy, leaving variable of the next minimum ratio also preserves the feasibility of the existing feasible variables.

If we concentrate in the value of $1b_{M_2}$ throughout the iterations, its value is strictly decreased for non-degenerate pivots and remained unchanged for degenerate pivots. So, finiteness of total number of bases in every LPP proves finiteness of our method for a complete non-degenerate LPP.                                ∎

The following example shows the comparison between traditional and our artificial free approaches.

## 9. A brief comparative study:

*Maximize Z  = 3x₁ + 5x₂*

$Maximize\ Z = 3x_1 + 5x_2$

*Subject to  $x_1 \leq 4$ , $x_2 \geq 6$, $x_1 + x_2 \geq 8$,  $3x_1 + 2x_2 \geq 18$*

*5x₁ + 4x₂ ≥ 32, x₁≥0, x₂≥0.*

$5x_1 + 4x_2 \geq 32,\ x_1 \geq 0,\ x_2 \geq 0.$

For the comparison purpose we first completely solve the above problem by traditional phase 1 (without using the trick mentioned in section 3) and then by our phase 1.

We assume that the reader should know how to construct the following initial dictionary of the traditional phase 1 method.

|     |     | x1 | x2 | v1 | S2 | s3 | s4 | s5 | Ratio |
|-----|-----|----|----|----|----|----|----|----|-------|
| z'  | -64 | -9 | -8 | 1  | 1  | 1  | 1  | 1  |       |
| Z   | 0   | -3 | -5 | 0  | 0  | 0  | 0  | 0  |       |
| s1  | 4   | 1* | 0  | -1 | 0  | 0  | 0  | 0  | 4     |
| v2  | 6   | 0  | 1  | 0  | -1 | 0  | 0  | 0  | --    |
| v3  | 18  | 3  | 2  | 0  | 0  | -1 | 0  | 0  | 6     |
| v4  | 8   | 1  | 1  | 0  | 0  | 0  | -1 | 0  | 8     |
| v5  | 32  | 5  | 4  | 0  | 0  | 0  | 0  | -1 | 32/5  |

Here $v_1 \geq 0$, $v_2 \geq 0$, $v_3 \geq 0$, $v_4 \geq 0$ and $v_5 \geq 0$ are artificial variables and z' is the phase 1 objective. Since $v_1$ is artificial non-basic variable, that means permanently non-basic, so we can safely remove its column. After performing simplex pivots we get the following sequence of dictionaries,

|     |     | s1 | x2 | s2 | s3 | s4 | s5 | Ratio |
|-----|-----|----|----|----|----|----|----|-------|
| z'  | -28 | 9  | -8 | 1  | 1  | 1  | 1  |       |
| z   | 12  | 3  | -5 | 0  | 0  | 0  | 0  |       |
| x1  | 4   | 1  | 0  | 0  | 0  | 0  | 0  | --    |
| v2  | 6   | 0  | 1  | -1 | 0  | 0  | 0  | 6     |
| v3  | 6   | -3 | 2* | 0  | -1 | 0  | 0  | 3     |
| v4  | 4   | -1 | 1  | 0  | 0  | -1 | 0  | 4     |

| v5 | 12 | -5 | 4 | 0 | 0 | 0 | -1 | 3 |
|---|---|---|---|---|---|---|---|---|

|  |  | s1 | s2 | s3 | s4 | s5 | Ratio |
|---|---|---|---|---|---|---|---|
| z' | -4 | -3 | 1 | -3 | 1 | 1 |  |
| z | 27 | - 9/2 | 0 | -5/2 | 0 | 0 |  |
| x1 | 4 | 1 | 0 | 0 | 0 | 0 | 4 |
| v2 | 3 | 3/2 | -1 | 1/2 | 0 | 0 | 2 |
| x2 | 3 | -3/2 | 0 | - 1/2 | 0 | 0 | -- |
| v4 | 1 | ½ | 0 | 1/2 | -1 | 0 | 2 |
| v5 | 0 | 1* | 0 | 2 | 0 | -1 | 0 |

|  |  | s2 | s3 | s4 | s5 | Ratio |
|---|---|---|---|---|---|---|
| Z' | -4 | 1 | 3 | 1 | -2 |  |
| Z | 27 | 0 | 13/2 | 0 | -9/2 |  |
| X1 | 4 | 0 | -2 | 0 | 1 | 4 |
| V2 | 3 | -1 | -5/2 | 0 | 3/2 | 2 |
| X2 | 3 | 0 | 5/2 | 0 | -3/2 | -- |
| V4 | 1 | 0 | -1/2 | -1 | ½* | 2 |
| S1 | 0 | 0 | 2 | 0 | -1 | -- |

|  |  | s2 | s3 | s4 | Ratio |
|---|---|---|---|---|---|
| z' | 0 | 1 | 1 | -3 |  |
| Z | 36 | 0 | 2 | -9 |  |
| x1 | 2 | 0 | -1 | 2 | 1 |
| v2 | 0 | -1 | -1 | 3* | 0 |
| x2 | 6 | 0 | 1 | -3 | -- |
| s5 | 2 | 0 | -1 | -2 | -- |
| s1 | 2 | 0 | 1 | -2 | -- |

|  |  | s2 | s3 |
|---|---|---|---|
| z' | 0 | 0 | 0 |
| z | 36 | -3 | -1 |
| x1 | 2 | 2/3 | - 1/3 |
| s4 | 0 | - 1/3 | - 1/3 |
| x2 | 6 | -1 | 0 |
| s5 | 2 | - 2/3 | 5/3 |
| s1 | 2 | - 2/3 | 1/3 |

It can be clearly observed that the sequence of visited corner points, in the original variable space, is (0,0),(4,0),(4,3),(2,6) and number of iterations are 5. The number of iterations are greater than number of visited corner points because of stalling due to degeneracy at the points (4,3) and (2,6).

Now, let us solve the same problem by using our new method. We constructed the following dictionary as described in section 2, by taking slack and surplus variables as initial basic variables. The vector $W$ could be obtained by directly summing over the infeasible constraints.   Our method shell provide the following sequence of dictionaries,

|     | W = | x1 [ -9 | x2 -8 ] | Ratio |
|-----|-----|---------|---------|-------|
| Z   | 0   | -3      | -5      |       |
| w1  | 4   | 1*      | 0       | 4     |
| w2  | -6  | 0       | -1      | ---   |
| w3  | -18 | -3      | -2      | 6     |
| w4  | -8  | -1      | -1      | 8     |
| w5  | -32 | -5      | -4      | 32/5  |

|     | W = | w1 [ 9 | x2 -8 ] | Ratio |
|-----|-----|--------|---------|-------|
| Z   | 12  | 3      | -5      |       |
| x1  | 4   | 1      | 0       | ---   |
| w2  | -6  | 0      | -1      | 6     |
| w3  | -6  | 3      | -2*     | 3     |
| w4  | -4  | 1      | -1      | 4     |
| w5  | -12 | 5      | -4      | 3     |

|     | W = | w1<br>[ -2 | w3<br>-1 ] | Ratio |
|-----|-----|------|------|-------|
| Z   | 27  | -7/2 | -3/2 |       |
| x1  | 4   | 1    | 0    | 4     |
| w2  | -3  | -3/2 | - 1/2 | 2    |
| x2  | 3   | -3/2 | - 1/2 | ---   |
| w4  | -1  | -1/2 * | - 1/2 | 2   |
| w5  | 0   | -1   | -2   | 0     |

|     | W = | w4<br>[ 0 | w3<br>0 ] |
|-----|-----|-----|-----|
| Z   | 36  | -9  | 2   |
| x1  | 2   | 2   | -1  |
| w2  | 0   | -3  | 1   |
| x2  | 6   | -3  | 1   |
| w4  | 2   | -2  | 1   |
| w5  | 2   | -2  | -1  |

Now it is clear that our method has the same sequence of visited corner points, as traditional phase 1 does, that is (0,0),(4,0),(4,3),(2,6) and the number of iterations are 3. Which are fewer than the previous approach because new method avoids stalling when the traditional phase 1 stalls due to degeneracy in artificial variables.

## 10. The artificial free dual simplex phase 1: (The dual counter part)

### Problem 2:

Given a dictionary $D(B)$, obtain dual feasibility.

### Algorithm :

**Step 1:** Let $K$ be a maximal subset of $N$ such that $K = \{j : d_{0j} < 0 , j \in N\}$. If $K = \emptyset$ then

$D(B)$ is dual feasible. **Exit.**

**Step 2:** Compute vector $W'(B) \in \Re^B$ such that $W'(B)_i = \sum_{k \in K} d_{ik}$ .

**Step 3:** Let $L \subseteq B$ such that $L = \{i : W'(B)_i < 0, \ i \in B\}$. If $L = \emptyset$ then $D(B)$ is dual

infeasible. **Exit.**

**Step 4:** Choose $r \in L$ such that $W'(B)_r \leq W'(B)_l \ \forall \ l \in L$

   (Ties should be broken arbitrarily)

**Step 5:** Choose $m \in N$ such that

   $m = \text{argmax} \ \{d_{0j}/d_{rj} \mid (d_{0j} < 0 \ , \ d_{rj} > 0) \text{ or } (d_{0j} \geq 0 \ , \ d_{rj} < 0), j \in N \}$

**Step 6:** Make a pivot on $(r, m)$. $(\Rightarrow$ Set $B := B + \{m\} - \{r\}$, $N := N - \{m\} + \{r\}$ and

   update $D(B)$.)

**Step 7:** Go to **Step 1**.

The above method is not more than just a dual counter-part of the algorithm described in section 6.

## 11. Conclusion:

This paper proposed equivalent but artificial free approaches of simplex and dual simplex phase 1 processes. The algorithm basically is a broad simplification of usual phase 1 simplex method and provides an advance in class room teaching. The new approach works in original variable space and obviates the use of artificial variables and constraints from respective traditional methods by allowing the negative variables into the basis, so, for any beginner it is very convenient to understand it.

Moreover another advantage is that because the new approach does not use artificial variables, it also avoids stalling whenever traditional phase 1 stalls due to degeneracy in artificial variables and reduces the storage requirement for the dictionary. Since it visits the same sequence of corner points as the traditional phase 1 does, its worst case complexity is also same as of the simplex method.